\documentclass{article}
\usepackage{amssymb}

\begin{document}
\title{Hydrodynamics for totally asymmetric $k$-step exclusion
processes}
\author{
H. Guiol\thanks{IMECC, Universidade de Campinas, Caixa Postal
6065, CEP  13053-970, Campinas, SP, Brasil},
K.Ravishankar\thanks{Dep. of Mathematics and Computer Science,
State University of New York, College at New Paltz, New Paltz, NY,
12561, USA}, E. Saada\thanks{CNRS, UPRES-A 6085, Universit\'e de
   Rouen, Site Colbert, 76821 Mont Saint Aignan Cedex,
   France}\\
AMS 1991 {\sl subject classifications}. Primary 60K35;Secondary
82C22.\\ {\sl Key words and phrases.} Hydrodynamics, $k$-step
exclusion,\\ nonconvex or nonconcave flux function, contact
discontinuities}

\maketitle


%
%

\begin{abstract}
We describe the hydrodynamic behavior of the $k$-step exclusion
process. Since the flux appearing in the hydrodynamic equation for
this particle system  is neither convex nor concave, the set of
possible solutions include in addition to entropic shocks and
continuous solutions those with contact discontinuities. We finish
with a limit theorem for the tagged particle.
\end{abstract}

\newcommand{\carn}{\hfill\rule{0.25cm}{0.25cm}}

\newtheorem{theorem}{Theorem}[section]
\newtheorem{lemma}[theorem]{Lemma}
\newtheorem{proposition}[theorem]{Proposition}
\newtheorem{corollary}[theorem]{Corollary}
\newtheorem{conjecture}[theorem]{Conjecture}
\newtheorem{definition}[theorem]{Definition}
\newtheorem{remark}[theorem]{Remark}
\def\E{{\mathbb E}}
\def\P{{\mathbb P}}
\def\R{{\mathbb R}}
\def\Z{{\mathbb Z}}
\def\V{{\mathbb V}}
\def\N{{\mathbb N}}

\section{Introduction and Notation}

In his paper Liggett (1980) introduced a Feller non conservative
approximation of the long range exclusion process to study the
latter. A conservative version of this dynamics, called $k$-step
exclusion process was defined and studied in Guiol (1999). It is
described in the following way.

Let $k\in\N^*:=\{1,2,...\}$, ${\bf X}:=\{0,1\}^{\Z}$ be the state
space, and let $\{X_n\}_{n\in{\N}}$ be a Markov chain on $\Z$ with
transition matrix $p(.,.)$ and ${\bf P}^x(X_0=x)=1$. Under the
mild hypothesis $\sup _{y\in {\Z}}\sum_{x\in {\Z}}p(x,y)<+\infty$,
$L_k$, defined below, is an infinitesimal pregenerator: For all
cylinder function $f$,

\begin{equation}\label{generateur}
L_kf(\eta )=\sum_{\eta (x)=1,\eta (y)=0}q_k(x,y,\eta )\left[
f(\eta^{x,y})-f(\eta )\right],
\end{equation}
where $q_k(x,y,\eta )={\bf E}^x\left[ \prod_{i=1}^{\sigma _y-1}
\eta (X_i),\sigma _y\leq\sigma_x,\sigma _y\leq k\right] $ is the
intensity for moving from $x$ to $y$ on configuration $\eta $,
$\sigma_y=\inf\left\{ n\geq 1:X_n=y\right\} $ is the first (non
zero) arrival time to site $y$ of the chain starting at site $x$
and $\eta^{x,y}$ is configuration $\eta$ where the states of sites
$x$ and $y$ were exchanged.

In words if a particle at site $x$ wants to jump it may go to the
first empty site encountered before returning to site $x$
following the chain $X_n$ (starting at $x$) provided it takes less
than $k$ attempts; otherwise the movement is cancelled.

By Hille-Yosida's theorem, the closure of $L_k$ generates a
continuous Markov semi-group $S_k(t)$ on $C({\bf X})$ which
corresponds to the $k$-step exclusion process $(\eta_t)_{t\geq
0}$. Notice that when $k=1$, $(\eta_t)_{t\geq 0}$ is the simple
exclusion process. An important property of $k$-step exclusion is
that it is an {\sl attractive process}.

Let ${\cal I}_k$ be the set of invariant measures for
$(\eta_t)_{t\geq 0}$ and let ${\cal S}$ be the set of translation
invariant measures on ${\bf X}$. If $p(x,y)=p(0,y-x)$ for all
$x,y\in{\Z}$ and $p(.,.)$ is irreducible then
\[
({\cal I}_k\cap {\cal S})_e=\left\{ \nu _\alpha :\alpha \in \left[
0,1\right] \right\},
\]
where the index $e$ mean extremal and $\nu_\alpha$ is the
Bernoulli product measure with constant density $\alpha$, {\it
i.e.} the measure with marginal
\[
\nu_{\alpha}\{\eta\in{\bf X}:\eta(x)=1\}=\alpha.
\]

In this paper we prove conservation of local equilibrium for the
totally asymmetric process in the Riemann case {\it i.e.}:
$p(x,x+1)=1$ for all $x\in\Z$ and the initial distribution is a
product measure with densities $\lambda$ to the left of the origin
and $\rho$ to its right, we denote it by $\mu_{\lambda,\rho}$. The
derived equation involves a flux which is neither concave nor
convex and appears for the first time as a hydrodynamic limit of
an interacting particle system. Up to now the ``constructive''
proofs for hydrodynamics relied on the concavity of the flux, see
Andjel \& Vares (1987) or the papers by Sepp\"al\"ainen ({\it
e.g.} Sepp\"al\"ainen (1998)), whose key tool is the Lax-Hopf
formula. The entropy solution for the type of equation we consider
was first studied in Ballou (1970). Such solutions can have
entropy shocks as well as contact discontinuities. Our aim is to
take advantage of Ballou's result to deduce conservation of local
equilibrium also in a constructive way, in the spirit of Andjel \&
Vares (1987).  However we explain in the last section how to
derive hydrodynamics for general initial profiles, and nearest
neighbor dynamics.

In section 3 we prove a law of large numbers for a tagged particle
in a $k$-step exclusion process.

\section{The hydrodynamic equation}

\subsection{Heuristic derivation of the equation}

Since the process is (totally) asymmetric we take Euler scaling.
For every $r\in\R$ define
$\eta_t^{\varepsilon}(r):=\eta_{\varepsilon^{-1}t}([\varepsilon^{-1}r])$,
where $[\varepsilon^{-1}r]$ is the integer part of
$\varepsilon^{-1}r$, for $\varepsilon>0$.

Given a continuous function $u^0(x),x\in\R$ (initial density
profile), we define a family of Bernoulli product measures
$\{\nu_{u^0}^\varepsilon\}_{\varepsilon>0}$ on ${\bf X}$, by: For
all $x \in\Z$,
\[
\nu_{u^0}^\varepsilon \{\eta\in{\bf X}:\eta
(x)=1\}=u^0(\varepsilon x).
\]
We call $\{\nu_{u^0}^\varepsilon\}$ the family of measures
determined by the profile $u^0$. Let $u^\varepsilon(r,t):=\int
(S_k(\varepsilon^{-1}t)\eta_0^\varepsilon(r))
d\nu^\varepsilon_{u^0}(\eta_0)$. Then for all $r\in\R$,
$u^\varepsilon(r,0)$ converges to $u^0(r)$  when $\varepsilon$
goes to 0; applying the generator (\ref{generateur}) to
$\eta_t^\varepsilon(r)$ we have
\[
\frac{d}{dt}S_k(t\varepsilon^{-1})\left(\eta_0^{\varepsilon}(r)\right)=
\varepsilon^{-1}S_k(t\varepsilon^{-1})\left[
-\sum_{i=0}^{k-1}\prod_{j=0}^i \eta_0^\varepsilon
(r+j\varepsilon)\left[1-\eta_0^\varepsilon\left(r+(i+1)\varepsilon
\right)\right] \right.
\]
\[
\left. +\sum_{i=0}^{k-1}\prod_{j=0}^i
\eta_0^\varepsilon(r-j\varepsilon)\left[1-\eta_0^\varepsilon(r)\right]
\right].
\]
Assuming that local equilibrium is preserved (thus expectation of
products factor), and taking expectations with respect to
$\nu_{u^0}^\varepsilon$, we obtain
\[
\frac{\partial u^\varepsilon}{\partial t} (r,t)=
\varepsilon^{-1}\left[ -\sum_{i=0}^{k-1}\prod_{j=0}^i
u^\varepsilon (r+j\varepsilon,t)\left[1-u^\varepsilon
\left(r+(i+1)\varepsilon, t\right) \right] \right.
\]
\[
\left. +\sum_{i=0}^{k-1}\prod_{j=0}^i
u^\varepsilon(r-j\varepsilon,t)\left[1-u^\varepsilon (r,t)\right]
\right].
\]

If we now let $\varepsilon$ converge to 0,
$u(r,t):=\lim_{\varepsilon\to 0}u^\varepsilon(r,t)$ should satisfy

\begin{equation}\label{hydro}
\left\{
\begin{array}{l}
\displaystyle{\frac{\partial u}{\partial t}+\frac{\partial
G_k(u)}{\partial x}=0} \\
 \\
u(x,0)=u^0(x),\\
\end{array}\right.
\end{equation}
where $G_k$ represents the flux of particles:
\[
G_k(u)=\sum_{j=1}^k ju^j(1-u).
\]

This is a non standard form because $G_k$ is neither convex nor
concave, thus equation (\ref{hydro}) is no longer ``a genuinely
nonlinear conservation law'', using the language of Lax (1973). To
deal with this equation, we have to use an extended version of non
linear Cauchy problems treated by Ballou (1970).

\begin{remark}
\end{remark}
Let $k$ go to infinity and denote by $G_{\infty}$ the limiting
flux function:
\[
G_{\infty}(u)=\frac u{1-u}.
\]
That case corresponds to the totally asymmetric long range
exclusion process. The resulting equation is simpler because the
flux function $G_{\infty}$ is strictly convex. The hydrodynamics
in this case should follow from the arguments of Aldous \&
Diaconis (1995) for the Hammersley's process.\par

\smallskip
For notational simplicity, from now on we restrict ourselves to
the case $k=2$. However our arguments can be easily extended for
all $k$.

\subsection{Hydrodynamics in the Riemann case}

\subsubsection{Notation and result}

Our main theorem characterizes the hydrodynamic (Euler) limit of
the 2-step exclusion process at points of continuity, when the
family of initial measures is determined by a step function
profile. From the heuristic derivation we would expect that in the
hydrodynamic limit the density profile would satisfy equation
(\ref{hydro}).  We show that the limiting density profile at time
$t$ is the entropy solution of equation (\ref{hydro}) starting
with the initial value $u^0$, a step function profile.  We now
give a brief summary of results concerning the solution of
equation (\ref{hydro}), due to D.P. Ballou (1970), when $u^0$ is a
step function. This will motivate the formulation of the theorem
as well as some aspects of the proof.

Existence of weak solution to the Cauchy problem given by equation
(\ref{hydro}) with bounded measurable initial condition was proved
in Ballou (1970), under the assumptions:

1. $G_k\in C^2(\R)$.

2. $G_k''$ vanishes at a finite number of points and changes sign
at these points.

In order to obtain uniqueness further conditions are needed.  We
require our solutions to satisfy the:

\noindent {\bf Condition E:} (O.A. Ole\u{\i}nik)

{\sl Let $x(t)$ be any curve of discontinuity of the weak solution
$u(t,x)$, and let $v$ be any number lying between
$u^-:=u(t,x(t)-0)$ and $u^+:=u(t,x(t)+0)$. Then except possibly
for a finite number of $t$,}
\[
S[v;u^-]\geq S[u^+;u^-],
\]
\noindent where
\[
S[v;w]:=\frac{G_k(w)-G_k(v)}{w-v}.
\]

It is known (Ballou (1970)) that the following two conditions are
necessary and sufficient for a piecewise smooth function $u(x,t)$
to be a weak solution of equation (\ref{hydro}):

1. $u(x,t)$ solves equation (\ref{hydro}) at points of smoothness.

2.If $x(t)$ is a curve of discontinuity of the solution then the
Rankine-Hugoniot condition ({\sl i.e.} $d (x(t))/dt = S[
u^+;u^-]$) holds along $x(t)$.

\noindent Moreover condition E is sufficient to ensure the
uniqueness of piecewise smooth solutions, which are the entropy
solutions to the equation. Hereafter we only deal with the case $k
= 2$. We denote $G(u):=G_2(u)$.\par

\smallskip
If $G$ were convex (concave) only two types of solutions would be
possible.  We now describe these two types of solutions.

Let $u^0(x)=\lambda 1_{\{x<0\}}+\rho 1_{\{x\geq 0\}}$.

If $\lambda>\rho$ ($\rho>\lambda$), then the speed of
characteristics which start from $x\leq 0$ (given by $G'$) is
greater than speed of characteristics which start from $x>0$. If
the intersection of characteristics occurs along a curve $x(t)$,
then since
\[
S[u^+;u^-]=\frac{G(\lambda)-G(\rho)}{\lambda-\rho}
=S[\lambda;\rho]
\]
Rankine-Hugoniot condition will be satisfied if
$x'(t)=S[\lambda;\rho]$.  Thus
\[
u(x,t)=\left\{
\begin{array}{ll}
\lambda, & x\leq S[\lambda;\rho]t;\\ \rho, & x >
S[\lambda;\rho]t.\\
\end{array}
\right.
\]
\noindent is a weak solution. The convexity of $G$ implies that
condition E is satisfied across $x(t)$. Therefore $u(x,t)$ defined
above is the unique entropic solution in this case and will be
referred to as a shock in the sense of Lax (1973).

\noindent If $\lambda<\rho \quad(\rho<\lambda)$, then the
characteristics starting respectively from $x\leq 0$ and from
$x>0$ never meet. Moreover they never enter the space-time wedge
between lines $x=\lambda t$ and $x=\rho t$. We can choose values
in this region to obtain a continuous solution, the so-called {\sl
continuous solution with a rarefaction fan}: Let $h$ be the
inverse of $H:=G'$,
\[
u(x,t)=\left\{
\begin{array}{ll}
\lambda, & x\leq H(\lambda)t;\\ h(x/t), & H(\lambda)t< x \leq
H(\rho)t;\\ \rho, & H(\rho)t< x.\\
\end{array}
\right.
\]
\noindent It is possible to define piecewise smooth weak solutions
with a jump occurring in the wedge satisfying the Rankine-Hugoniot
condition. But the convexity of $G$ prevents such solutions to
satisfy condition E. Thus  the continuous solution with a
rarefaction fan is the unique entropic solution in this case.\par

\smallskip
For the 2-step exclusion process, the flux function $G$ is neither
concave nor convex.  Instead $G(u)=u+u^2- 2 u^3$ is convex for
$u<1/6$ and concave for $u>1/6$.  In this case in addition to the
shock and continuous solution with a rarefaction fan it is
possible to have solutions for which the curve of discontinuities
never enters the region of intersecting characteristics. The
quotation in boldface is the original number of lemmas, prop... in
Ballou (1970), but the notation refers to 2-step exclusion:

\begin{definition}{\bf [B def2.1]}\label{Bd2.1}

For any $u< 1/6$, define $u^*:=u^*(u)$ as
\[
u^*=\sup\{\eta>u: S[u;\eta]>S[v;u]\ \forall v\in(u,\eta)\}.
\]

For any $u> 1/6$, define $u_*:=u_*(u)$ as
\[
u_*=\inf\{\eta<u: S[u;\eta]>S[v;u]\ \forall v\in(\eta,u)\}.
\]
\end{definition}

In other words, for $u< 1/6$, if we consider the upper convex
envelope $G^c$ of $G$ on $(u,+\infty)$, then $u^*$ is the first
point where $G^c$ coincides with $G$. In the same way when $u>
1/6$, $u_*$ is the first point where the lower convex envelope
$G_c$ of $G$ on $(-\infty,u)$ coincides with $G$. For $\eta< 1/6$,
$\eta^*=(1-2\eta)/4$, and for $\eta> 1/6$, $\eta_*=(1-2\eta)/4$.

Let $h_1$ and $h_2$ be the inverses of $H$ respectively restricted
to $(-\infty,1/6)$ and to $(1/6,+\infty)$, {\it i.e.}
$h_1(x)=(1/6)(1-\sqrt{7-6x})$ and $h_2(x)=(1/6)(1+\sqrt{7-6x})$
for $x\in(-\infty,7/6)$.

The following lemmas are taken from Ballou (1970).

\begin{lemma}{\bf [B lem2.2]}\label{Bl2.2}
Let $\eta<1/6$ be given, and suppose that $\eta^*<\infty$. Then
$S[\eta;\eta^*]=H(\eta^*)$.
\end{lemma}

\begin{lemma}{\bf [B lem2.4]}\label{Bl2.4}
Let $\eta< 1/6$ be given, and suppose that $\eta^*<\infty$. Then
$\eta^*$ is the only zero of $S[u;\eta]-H(u)$, $u>\eta$.
\end{lemma}

If $\lambda<\rho<1/6$, the relevant part of the flux function is
convex and the unique entropic weak solution is the continuous
solution with a rarefaction fan.

If $\rho <\lambda<\rho^*$ $(\rho <1/6)$, then $H(\eta)>H(\rho)$ if
$\rho<\eta\leq 1/6$, and $H(\eta)>H(\rho^*)>H(\rho)$ if
$1/6<\eta<\rho^*$ since $H$ is decreasing in this region. Thus
$H(\lambda)>H(\rho)$, which implies an intersection of
characteristics: The unique entropic weak solution is the shock.

Let $\rho<\rho^*< \lambda$ $(\rho<1/6)$: Lemma \ref{Bl2.4} applied
to $\rho$ suggests that a jump from $\rho^*$ to $\rho$ along the
line $x=H(\rho^*)t$ will satisfy the Rankine-Hugoniot condition.
Since $\rho^*$ is specially defined for this, a solution with such
a jump will also satisfy condition E. Therefore if we can
construct a solution with the jump described above it will be the
unique entropic weak solution in this case. Notice that since
$H(\lambda)<H(\rho^*)$, no characteristics intersect along the
line of discontinuity $x=H(\rho^*)t$.  We call this case {\sl
contact discontinuity}, following Ballou.  The solution is defined
by
\[
u(x,t)=\left\{
\begin{array}{ll}
\lambda, & x\leq H(\lambda)t;\\ h_2(x/t), & H(\lambda)t<x\leq
H(\rho^*)t;\\ \rho, & H(\rho^*)t< x.\\
\end{array}
\right.
\]
Corresponding cases on the concave side of $G$ are treated
similarly.\par
\smallskip
Let $\tau$ denote the shift operator. We are able now to state our
result.

\begin{theorem}\label{limhydro}
Let $v\in\R$, $\lambda,\rho\neq 1/6$, and $\mu_{\lambda,\rho}$ the
Bernoulli product measure on $\Z$ with densities $\lambda$ for
$x\leq 0$ and $\rho$ for $x>0$. Then
\[
\lim_{t\to\infty} \mu_{\lambda,\rho}\tau_{[vt]}S_2(t)=\nu_{u(v,1)}
\]
at every continuity point of $u(.,1)$, where $\nu_{u(v,1)}$
denotes the product measure with density $u(v,1)$ defined by:

Case 1. $\lambda<\rho< 1/6$: continuous solution, with a
rarefaction fan
\[
u(x,1)=\left\{
\begin{array}{ll}
\lambda, & x\leq H(\lambda);\\ h_1(x), & H(\lambda)< x\leq
H(\rho);\\ \rho, & H(\rho)< x.\\
\end{array}
\right.
\]

Case 2. $\rho<\lambda<\rho^*$, ($\rho< 1/6$): entropy shock
\[
u(x,1)=\left\{
\begin{array}{ll}
\lambda, & x\leq S[\lambda;\rho];\\ \rho, & x> S[\lambda;\rho].\\
\end{array}
\right.
\]

Case 3. $\rho<\rho^*<\lambda$, ($\rho< 1/6$): contact
discontinuity
\[
u(x,1)=\left\{
\begin{array}{ll}
\lambda, & x\leq H(\lambda);\\ h_2(x), & H(\lambda)< x\leq
H(\rho^*);\\ \rho, & H(\rho^*)< x.\\
\end{array}
\right.
\]

Case 4. $1/6<\rho<\lambda$: continuous solution, with a
rarefaction fan
\[
u(x,1)=\left\{
\begin{array}{ll}
\lambda, & x\leq H(\lambda);\\ h_2(x), & H(\lambda)< x\leq
H(\rho);\\ \rho, & H(\rho)< x.\\
\end{array}
\right.
\]

Case 5. $\rho>\lambda>\rho_*$, ($\rho>1/6$): entropy shock
\[
u(x,1)=\left\{
\begin{array}{ll}
\lambda, & x\leq S[\lambda;\rho];\\ \rho, & x> S[\lambda;\rho].\\
\end{array}
\right.
\]

Case 6. $\rho>\rho_*>\lambda$, ($\rho> 1/6$): contact
discontinuity
\[
u(x,1)=\left\{
\begin{array}{ll}
\lambda, & x\leq H(\lambda);\\ h_1(x), & H(\lambda)< x\leq
H(\rho_*);\\ \rho, & H(\rho_*)< x.\\
\end{array}
\right.
\]

\end{theorem}

\begin{remark}
\end{remark}
For any $k\geq 2$ the profiles will be of the same kind, because
$G_k$ has only one inflection point between 0 and 1 and is first
convex then concave.

\begin{remark}
\end{remark}
Comparing with hydrodynamics of simple exclusion we observe that
$k$-step exclusion ($k\geq 2$) has not only a stable increasing
shock (Case 5) and a decreasing continuous solution (Case 1) but
also a stable decreasing shock (Case 2), an increasing continuous
solution (Case 4) and two contact discontinuities (Cases 3 and 6).

\subsubsection{Proof of Theorem \ref{limhydro}}

It follows the scheme introduced in Andjel \& Vares (1987), where
the authors obtained the hydrodynamic limit for the
one-dimensional zero-range process in the Riemann case, i.e. the
hydrodynamic equation
\[
\left\{
\begin{array}{ll}
\displaystyle{\frac{\partial u}{\partial t}+\frac{\partial
\phi(u)}{\partial
    x}=0}\\
 \\
 u(x,0)=u^0(x)=\lambda 1_{\{x<0\}}+\rho 1_{\{x\geq 0\}}\\
\end{array}
\right.
\]
was derived. There $\phi$, the mean flux of particles through the
origin, was a concave function. Therefore, their proof used both
the monotonicity of the process, still valid here, and the
concavity of the flux, that we have to replace by an {\it ad hoc}
use of the properties of the solution of (\ref{hydro}).

Informally speaking, they first showed that a weak Ces\'aro limit
of (the measure of) the process is an invariant and translation
invariant measure. Then they showed that the (Ces\'aro) limiting
density inside a macroscopic box is equal to the difference of the
edge values of a flux function. These propositions were based on
monotonicity, and on the characterization of invariant and
translation invariant measures (both valid for k-step as
well),thus we can quote them (with appropriate notation for the
$k$-step), and take them for granted.

\begin{lemma}{\bf [AV 3.1]}
Let $\mu$ be a probability measure on $\{0,1\}^{\Z}$ such that

(a) $\nu_\rho\leq\mu\leq\nu_\lambda$ for some
$0\leq\rho<\lambda\leq 1$,\ \ (b) either $\mu\tau_1\leq\mu$ or
$\mu\tau_1\geq\mu$.

\noindent Then any sequence $T_n\to\infty$ has a subsequence
$T_{n_k}$ for which there exists $D$ dense (countable) subset of
$\R$ such that for each $v\in D$,
\[
\lim_{k\to\infty}\frac 1 {T_{n_k}} \int_0^{T_{n_k}}
\mu\tau_{[vt]}S_2(t)dt=\mu_v
\]
for some $\mu_v\in{\cal I}_2\cap {\cal S}$.
\end{lemma}

\begin{lemma}{\bf [AV 3.2]}
For $v\in D$, we can write
$\mu_v=\int\nu_{\alpha}\gamma_v(d\alpha)$, where $\gamma_v$ is a
probability on $[\rho,\lambda]$. Also, if $u<v$ are in $D$,
\begin{equation}\label{3.4}
\lim_{k\to\infty} \mu S_2({T_{n_k}}) \left(\frac 1
{T_{n_k}}\sum_{[uT_{n_k}]}^{[vT_{n_k}]}\eta(x)\right)= F(v)-F(u)
\end{equation}
with, for $w\in D$,
$F(w)=\int[w\alpha-G(\alpha)]\gamma_w(d\alpha).$
\end{lemma}

The difficult part is then to prove that $\gamma_v$ is in fact the
Dirac measure concentrated on $u(x,1)$. They did it in Lemma [AV
3.3] and Theorem [AV 2.10] using the concavity of their flux
function.

For $k$-step exclusion, at this point we will have to look
separately at the six cases given in the theorem.

To conclude, Andjel \& Vares had to prove that the Ces\'aro limit
implies a weak limit, through the following propositions, based on
monotonicity. We will use these results also without proof.

\begin{proposition}{\bf [AV 3.4]}\label{AV3.4}
Let $\mu=\mu_{\lambda,\rho}$. If
\[
\mu_v=\left\{
\begin{array}{ll}
\nu_\lambda, & {\it if\,}v\in D,\, v< S[\lambda;\rho]\\ \nu_\rho,
& {\it if\,}v\in D,\, v> S[\lambda;\rho]\\
\end{array}
\right.
\]
then
\[
\lim_{T\to\infty}\frac 1 T \int_0^T \mu\tau_{[vt]}S_2(t)dt=\mu_v.
\]
\end{proposition}

\begin{proposition}{\bf [AV 3.5]}\label{AV3.5}
If $\mu$ satisfies

(a) $\mu\leq\nu_\lambda$, (b) $\mu\tau_1\geq\mu$, (c) there exists
$v_0$ finite so that
\[
\lim_{T\to\infty}\frac 1 T \int_0^T
\mu\tau_{[vt]}S_2(t)dt=\nu_\lambda
\]
for all $v>v_0$. Then
\[
\lim_{t\to\infty}\mu\tau_{[vt]}S_2(t)=\nu_\lambda \quad {\rm
for\,all\,\,}v>v_0.
\]
\end{proposition}

\noindent {\bf Proof of Theorem \ref{limhydro} in Cases 2 and 3.}
\noindent In case 1, the proof is not different from the one given
in Andjel \& Vares (1987) since $G$ is convex in the relevant
region; thus we omit it. In case 2, $G$ is not convex in the
relevant region, therefore we supply a proof, though it is quite
close to the original one. Case 3 uses ideas from cases 1 and 2
and introduces some new ideas to deal with complications arising
from the non-convexity of $G$.  Cases 4-6 are symmetric to cases
1-3 in the sense that the roles played by convexity and concavity
are exchanged. \par

\smallskip
The proof has 3 steps. The main ingredients are monotonicity, and
inequalities relying on the properties of $S[.;.],\,G,\,H$ given
before.

\noindent {\it First step.}\par Using monotonicity of the 2-step
exclusion, we can proceed as in the beginning of the proof of
Lemma [AV3.3], and get two finite values $\underbar v$ and $\bar
v$ so that: If $v\in D$ and $v>\bar v$, then
$\gamma_v=\delta_\rho$, while $\gamma_v=\delta_\lambda$ if
$v<\underbar v$.

\smallskip
\noindent {\it Second step, preliminary.}\par Let $u<v$, both in
$D$. Attractiveness of the process (since
$\nu_\rho\leq\mu_{\lambda,\rho}\leq\nu_\lambda$) and (\ref{3.4})
imply
\begin{equation}\label{3.12}
\begin{array}{c}
\displaystyle{(v-u)\rho\leq v\int\alpha\gamma_v(d\alpha)-\int
G(\alpha)\gamma_v(d\alpha) -u\int\alpha\gamma_u(d\alpha)+\int
G(\alpha)\gamma_u(d\alpha)}\\ \displaystyle{\leq (v-u)\lambda}\\
\end{array}
\end{equation}

(i) Taking $u<\underbar v$, the first step gives
$\gamma_u=\delta_\lambda$, so the second inequality of
(\ref{3.12}) is simplified in
\[
v\int\alpha\gamma_v(d\alpha)-\int G(\alpha)\gamma_v(d\alpha)
-u\lambda+G(\lambda) \leq (v-u)\lambda
\]
which can be written
\begin{equation}\label{*1}
\int_{[\rho,\lambda]}(G(\lambda)-G(\alpha))\gamma_v(d\alpha) \leq
v\int_{[\rho,\lambda]}(\lambda-\alpha)\gamma_v(d\alpha)
\end{equation}

(ii)  Similarly, for $\bar v<v$, $\gamma_v=\delta_\rho$, and the
first inequality of (\ref{3.12}) reads
\begin{equation}\label{*2}
u\int_{[\rho,\lambda]}(\alpha-\rho)\gamma_u(d\alpha)\leq
\int_{[\rho,\lambda]}(G(\alpha)-G(\rho))\gamma_u(d\alpha)
\end{equation}

\noindent {\bf Proof for Case 2:} $\rho<\lambda
<\rho^*,\,\rho<1/6.$\par The definition of $\rho^*$ implies that
for every $\alpha\in (\rho,\lambda)$
\[
\frac{G(\lambda) - G(\alpha)}{\lambda - \alpha} \geq \frac
{G(\lambda) - G(\rho)}{\lambda - \rho} = S[\lambda;\rho],
\]
so that inequality (\ref{*1}) for $\underbar v\leq
v<S[\lambda;\rho]$ ($v\in D$) yields
\[
v\int_{[\rho,\lambda]}(\lambda -\alpha) \gamma_v (d\alpha) \geq
S[\lambda;\rho]\int_{[\rho,\lambda]}(\lambda-\alpha)\gamma_v
(\alpha).
\]
Since $v< S[\lambda;\rho]$ we conclude that $\gamma_v =
\delta_\lambda$.

Starting from inequality (\ref{*2}) with $\bar v\geq u
>S[\lambda;\rho]$ ($u\in D$) and proceeding in a similar manner we
can show that $\gamma_u = \delta_\rho$.\par

\smallskip
\noindent {\bf Proof for case 3:} $\rho<\rho^*<\lambda,\,\rho<
1/6$.\par \noindent {\it Second step, part 1.}\par Let $u<v$, both
in $D$, $u<\underbar v$, and $v<H(\lambda)$. On
$[\rho^*,\lambda)$, $G$ is concave, thus
 for every $\alpha\in[\rho^*,\lambda]$
\begin{equation}\label{conc1}
H(\lambda)<S[\alpha;\lambda] \le H(\rho^*).
\end{equation}
If $\alpha \in [\rho,\rho^*)$ then from the definition of $\rho^*$
it follows that
\begin{equation}\label{conc2}
S[\alpha;\rho^*] \geq H(\rho^*).
\end{equation}
We decompose $[\rho,\lambda]=[\rho,\rho^*]\cup(\rho^*,\lambda]$,
and (\ref{*1}) becomes
\begin{equation}\label{*6}
\begin{array}{c}
\displaystyle{v\int_{[\rho,\rho^*]}(\lambda-\rho^*)\gamma_v(d\alpha)
+v\int_{[\rho,\rho^*]}(\rho^*-\alpha)\gamma_v(d\alpha)
+v\int_{(\rho^*,\lambda]}(\lambda-\alpha)\gamma_v(d\alpha) \geq}
\\
\noindent \displaystyle{\int_{[\rho,\rho^*]}
(G(\lambda)-G(\alpha))\gamma_v(d\alpha)+
\int_{(\rho^*,\lambda]}(G(\lambda)-G(\alpha))\gamma_v(d\alpha)}.\\
\end{array}
\end{equation}
By (\ref{conc2})
\[
\int_{[\rho,\rho^*]}(G(\lambda)-G(\alpha))\gamma_v(d\alpha)=
\int_{[\rho,\rho^*]}
\left[(G(\lambda)-G(\rho^*)+G(\rho^*)-G(\alpha)\right]\gamma_v(d\alpha)
\]
\[
\geq\int_{[\rho,\rho^*]} \left[(\lambda-\rho^*)S[\rho^*;\lambda]
+H(\rho^*) (\rho^*-\alpha)\right]\gamma_v(d\alpha)
\]
and by (\ref{conc1})
\[
\int_{(\rho^*,\lambda]}(G(\lambda)-G(\alpha))\gamma_v(d\alpha)\geq
H(\lambda)\int_{(\rho^*,\lambda]}(\lambda-\alpha)\gamma_v(d\alpha).
\]
Those two inequalities together with (\ref{*6}) give
\[
(v-S[\rho^*;\lambda])(\lambda-\rho^*)
\int_{[\rho,\rho^*]}\gamma_v(d\alpha)+ (v-H(\rho^*))
\int_{[\rho,\rho^*]}(\rho^*-\alpha)\gamma_v(d\alpha)
\]
\[
+(v-H(\lambda))
\int_{(\rho^*,\lambda]}(\lambda-\alpha)\gamma_v(d\alpha)\geq 0.
\]
Since $v<H(\lambda)$, by (\ref{conc1}), the only possibility is
$\gamma_v = \delta_{\lambda}$.

\bigskip
\noindent {\it Second step, part 2.}\par Let $u,v\in D$, such that
$v>\bar v$, and $H(\rho^*)<u<v$. We decompose each integral of
(\ref{*2}) on the two intervals $[\rho,\rho^*]$ and
$(\rho^*,\lambda]$.
\[
\int_{(\rho^*,\lambda]}(G(\alpha)-G(\rho))\gamma_u(d\alpha)
\]
\[
=\int_{(\rho^*,\lambda]}(G(\alpha)-G(\rho^*))\gamma_u(d\alpha)+
\int_{(\rho^*,\lambda]}(G(\rho^*)-G(\rho))\gamma_u(d\alpha)
\]
By definition of $\rho^*$ and by Lemma \ref{Bl2.2}, for all
$\alpha \in [\rho,\rho^*]$,
\[
G(\alpha)-G(\rho) \le H(\rho^*)(\alpha-\rho),
\]

and $G$ being strictly concave on $(\rho^*,\lambda]$, we have for
all $\alpha \in [\rho^* ,\lambda]$
\[
G(\alpha)-G(\rho^*) \le H(\rho^*)(\alpha-\rho^*).
\]
Thus by (\ref{*2})
\[
\begin{array}{c}
\displaystyle{u\int_{[\rho,\rho^*]}(\alpha-\rho)\gamma_u(d\alpha)+
u\int_{(\rho^*,\lambda]}(\alpha-\rho^*)\gamma_u(d\alpha)+
u\int_{(\rho^*,\lambda]}(\rho^*-\rho)\gamma_u(d\alpha)}\\
\displaystyle{\le \int_{[\rho,\rho^*]}
H(\rho^*)(\alpha-\rho)\gamma_u(d\alpha)+
\int_{(\rho^*,\lambda]}H(\rho^*)(\alpha-\rho^*)\gamma_u(d\alpha)}\\
\displaystyle{+\int_{(\rho^*,\lambda]}
H(\rho^*)(\rho^*-\rho)\gamma_u(d\alpha)}\\
\end{array}
\]
which can be written
\[
\begin{array}{c}
\displaystyle{(u-H(\rho^*))\left[\int_{[\rho,\rho^*]}(\alpha-\rho)
\gamma_u(d\alpha)+\int_{(\rho^*,\lambda]}
(\alpha-\rho^*)\gamma_u(d\alpha) \right.}\\
\displaystyle{\left.+\int_{(\rho^*,\lambda]} (\rho^*-\rho)
\gamma_u(d\alpha)\right] \le 0}.\\
\end{array}
\]
Since $u-H(\rho^*)>0$, this implies
\[
\int_{[\rho,\lambda]}(\alpha-\rho)\gamma_u(d\alpha)= 0
\]

\noindent Therefore we conclude $\gamma_u=\delta_\rho$.

\bigskip
\noindent {\it Second step, conclusion.}\par Using the two
preceding parts, attractivity, Propositions \ref{AV3.4} and
\ref{AV3.5}, we conclude in case 3
\begin{equation}\label{4.2.2}
\lim_{t\to\infty}\mu_{\lambda,\rho}\tau_{[vt]}S_2(t)=\left\{
\begin{array}{ll}
\nu_\lambda, & {\it if\,\,} v< H(\lambda)\\ \nu_\rho, & {\it
if\,\,} v> H(\rho^*)\\
\end{array}
\right.
\end{equation}
\noindent and in case 2
\begin{equation}\label{4.2.2a}
\lim_{t\to\infty}\mu_{\lambda,\rho}\tau_{[vt]}S_2(t)=\left\{
\begin{array}{ll}
\nu_\lambda, & {\it if\,\,} v< S[\lambda;\rho]\\ \nu_\rho, & {\it
if\,\,} v> S[\lambda;\rho]\\
\end{array}
\right.
\end{equation}

\bigskip
\noindent {\it Third step, first part.}\par Let $u_1,v,v_1\in D$,
$u_1<H(\lambda)<v<H(\rho^*)<v_1$, such that $u_1$ and $v_1$ belong
to an interval where $G$ is concave. By attractivity,
\begin{equation}\label{aux3}
\limsup_{t\to\infty}\mu_{\lambda,\rho} S_2(t) \left(\frac 1
t\sum_{[u_1 t]}^{[ v t]}\eta(x)\right) \leq(v-u_1)\lambda
\end{equation}
The second step and (\ref{3.4}) imply
\[
\lim_{t\to\infty}\mu_{\lambda,\rho} S_2(t) \left(\frac 1
t\sum_{[u_1 t]}^{[v_1 t]}\eta(x)\right)
=v_1\rho-G(\rho)-u_1\lambda+G(\lambda)
\]
which, combined with (\ref{aux3}) gives
\begin{equation}\label{4.5.2}
\liminf_{t\to\infty}\mu_{\lambda,\rho} S_2(t) \left(\frac 1
t\sum_{[v t]}^{[v_1 t]}\eta(x)\right) \geq
v_1\rho-G(\rho)+G(\lambda)-v\lambda
\end{equation}
Since $H(\lambda)<v<H(\rho^*)$, there exists some $\theta$,
$\rho^*<\theta<\lambda$, with $v=H(\theta)$, so that
$\theta=h_2(v)$. Let $\theta<\theta'<\lambda$, we now apply
(\ref{4.5.2}) to $\mu_{\theta',\rho}$, and we use attractivity
through $\mu_{\theta',\rho}\leq\mu_{\lambda,\rho}$ to get
\begin{equation}\label{aux4}
\begin{array}{c}
\displaystyle{\liminf_{t\to\infty}\mu_{\lambda,\rho} S_2(t)
\left(\frac 1 t\sum_{[v t]}^{[v_1 t]}\eta(x)\right)
\geq\liminf_{t\to\infty}\mu_{\theta',\rho} S_2(t) \left(\frac 1
t\sum_{[v t]}^{[v_1 t]}\eta(x)\right)}\\
 \\
\displaystyle{\geq v_1\rho-G(\rho)+G(\theta')-v\theta'}\\
\end{array}
\end{equation}
By (\ref{3.4}) again,
\[
\liminf_{t\to\infty}\mu_{\lambda,\rho} S_2(t) \left(\frac 1
t\sum_{[v t]}^{[v_1 t]}\eta(x)\right) \leq v_1\rho-G(\rho)
-v\int_{[\rho,\lambda]}\alpha\gamma_v(d\alpha)
+\int_{[\rho,\lambda]} G(\alpha)\gamma_v(d\alpha)
\]
which, together with (\ref{aux4}), if we make $\theta'\to\theta$,
gives
\begin{equation}\label{4.7.2}
G(\theta)- v\theta\leq
\int_{[\rho,\lambda]}(G(\alpha)-v\alpha)\gamma_v(d\alpha)
\end{equation}
Since $G(\alpha)$ is convex when $\alpha\in (\rho,1/6)$ and
concave when $\alpha\in (1/6,\lambda)$, for all $v$,
$G(\alpha)-v\alpha$  has at most two critical points determined by
the condition $G'(\alpha)=v$. One satisfies $1/6<\alpha<\lambda$
and is a local maximum, the other (when it exists) is a local
minimum with $\rho<\alpha<1/6$. We want to conclude that the local
maximum is a global maximum when $\rho<\alpha<\lambda$. We know
from Lemma (\ref{Bl2.2}) that
\[
G(\rho^*) - G(\rho) = (\rho^* -\rho) H(\rho^*)
\]
which implies, since $v=H(\theta)<H(\rho^*)$,
\[
G(\rho^*) - v \rho^*  > G(\rho) - v \rho
\]
Because $G(\alpha)-v\alpha$ is increasing in $(\rho^*,\theta)$
(recall that $H(\theta) = v$) we have: For all $\alpha\in
(\rho^*,\theta)$
\[
G(\alpha)-v\alpha>G(\rho^*)-v\rho^*>G(\rho)-v\rho
\]
We conclude that $\theta$ is a global maximum:
\[
\max_{\rho\leq\alpha\leq\lambda}[G(\alpha)-v\alpha]=G(\theta)-v\theta
\]
thus $\gamma_v=\delta_\theta=\delta_{h_2(v)}$.

\bigskip
\noindent {\it Third step, second part.}\par This part follows
closely the argument in Andjel \& Vares (1987), but we detail it
for the sake of completeness. Since the measures $\frac 1 T
\int_0^T \mu_{\lambda,\rho}\tau_{[vt]}S_2(t)\ dt$ depend
monotonically on $v$ and form a relatively compact set, we have
for all $v$,
\[
\lim_{T\to\infty} \frac 1 T \int_0^T
\mu_{\lambda,\rho}\tau_{[vt]}S_2(t)\ dt=\nu_{h_2(v)}.
\]
It remains to prove that
\[
\lim_{t\to\infty}\mu_{\lambda,\rho}\tau_{[vt]}S_2(t)=\nu_{h_2(v)}
\]
when $H(\lambda)<v<H(\rho^*)$ (by continuity of $h_2$, this result
will also be valid at $H(\lambda)$). For this, let $\tilde\mu_v$
be a weak limit of $\mu_{\lambda,\rho}\tau_{[vt]}S_2(t)$. It is
enough to show
\[
({\it a})\,\tilde\mu_v\geq\nu_{h_2(v)},\qquad ({\it
b})\,\tilde\mu_v(\eta(0))=h_2(v).
\]

(a) Let $\theta=h_1(v)$, with $H(\lambda)<v<H(\rho^*)$,
$v=H(\theta)$. Let $\rho^*<\tilde\theta<\theta<\lambda$, then
$v<H(\tilde\theta)$ and
$\mu_{\tilde\theta,\rho}\tau_{[vt]}S_2(t)\leq
\mu_{\lambda,\rho}\tau_{[vt]}S_2(t)$. Added to (\ref{4.2.2}), this
yields
\[
\lim_{t\to\infty}
\mu_{\tilde\theta,\rho}\tau_{[vt]}S_2(t)=\nu_{\tilde\theta}\leq
\lim_{t\to\infty}\mu_{\lambda,\rho}\tau_{[vt]}S_2(t)=\tilde\mu_v.
\]
Hence, by continuity, if $\tilde\theta$ converges to $\theta$,
\[
\tilde\mu_v\geq\nu_\theta=\nu_{h_2(v)}.
\]

(b) Let $u_1<H(\lambda)<H(\rho^*)<v_1$. By the definition of
$u(x,t)$ in that case,
\[
\int_{u_1}^{v_1} u(x,1)\ dx = \int_{u_1}^{H(\lambda)}\lambda\ dx+
\int_{H(\lambda)}^{H(\rho^*)}h_2(x)\ dx+
\int_{H(\rho^*)}^{v_1}\rho\ dx
\]
\[
= \lambda(H(\lambda)-u_1)+ \int_\lambda^{\rho^*} \theta
H'(\theta)\ d\theta +\rho(v_1-H(\rho^*))
\]
\[
=v_1\rho-u_1\lambda -G(\rho)+G(\lambda)
\]
where we have integrated by parts the integral of the second line
(derived by a change of variables); therefore, by (\ref{3.4}),
\begin{equation}\label{4.12.2}
\lim_{t\to\infty}\mu_{\lambda,\rho} S_2(t) \left(\frac 1
t\sum_{[u_1 t]}^{[v_1 t]}\eta(x)\right)= \int_{u_1}^{v_1} u(x,1)\
dx
\end{equation}
and (b) follows from (\ref{4.12.2}) and (a), with the same
argument as in \cite{andkip} p.332 (proof of Theorem 3.2), based
on monotonicity.$\carn$

\section{Asymptotic behavior of a tagged particle}


We introduce here an interpretation of the $k$-step exclusion
dynamics valid in the totally asymmetric case. Up to now we
considered that a particle might jump from $x$ to the first empty
site in $\{x+1,...,x+k\}$. If we want to leave the particles
ordered we could equally say that the particle at $x$ pushes the
``pack'' of ($\leq k$) neighboring particles in front of it, each
one moving of one unit to the right. In other words, if there is
no particle at site $x+1$ then the particle at site $x$ goes to
site $x+1$; if there is one particle at site $x+1$ and no particle
at site $x+2$ then the particle at $x$ pushes the particle of site
$x+1$ to site $x+2$ and occupies site $x+1$; and so on... Then the
generator reads
\begin{equation}\label{pushing}
L_k f(\eta)=\sum_{x\in{\bf Z}}\sum_{i=0}^{k-1}
\prod_{j=0}^i\eta(x+j)\left(1-\eta(x+i+1)\right)
\left[f(\eta^{x,x+1,...,x+i+1})-f(\eta)\right]
\end{equation}
where
\[
\eta^{x_1,x_2,...,x_l}(u)=\left\{
\begin{array}{ll}
\eta(x_l)&\mbox{ if }u=x_1,\\ \eta(x_{i-1})&\mbox{ if }u=x_i,\
i=2,...,l,\\ \eta(u)&\mbox{ otherwise}.
\end{array}
\right.
\]

It is easy to see, comparing (\ref{pushing}) to
(\ref{generateur}), that they do correspond in this setting,
because we do not label the particles. The interest of
(\ref{pushing}) w.r.t. (\ref{generateur}) is that it keeps track
of the particles' order. We will need this interpretation in the
next section.

Similarly, we define a Tagged ``Pushing'' Particle, and the
generator of the $k$-step exclusion process as seen from this
Tagged Pushing Particle is

\[
\begin{array}{l}
\widetilde{L}_k(\eta)=\\ \displaystyle{\sum_{i=0}^{k-1}\sum_{x\neq
  0,-1,...,-(i+1)}\prod_{j=0}^i\eta(x+j)\left(1-\eta(x+i+1)\right)
\left[ f(\eta^{x,x+1,...,x+i+1})-f(\eta) \right]}\\
\displaystyle{+\sum_{n=1}^k \prod_{m=1}^{n-1}
\eta(m)(1-\eta(n))[f(\tau_1\eta^{0,1,...,n})-f(\eta)]}\\
\displaystyle{+\sum_{n=1}^{k-1}\sum_{l=1}^{k-n}\prod_{m=-l}^{-1}
\eta(m)\prod_{i=1}^{m-1} \eta(i)(1-\eta(n))
\left[f(\tau_1\eta^{-l,...,0,...,n})-f(\eta)\right].}\\
\end{array}
\]
\noindent To be clearer, let us write and comment it for $k=2$.

\begin{equation}\label{un}
\widetilde{L}_2(\eta)=\sum_{x\neq 0,-1}\eta(x)(1-\eta(x+1)) \left[
f(\eta^{x,x+1})-f(\eta) \right]
\end{equation}
\begin{equation}\label{deux}
+(1-\eta(1))\left[ f(\tau_1\eta^{0,1})-f(\eta) \right]
\end{equation}
\begin{equation}\label{trois}
+\sum_{x\neq 0,-1,-2}\eta(x)\eta(x+1)(1-\eta(x+2))
\left[f(\eta^{x,x+1,x+2})-f(\eta)\right]
\end{equation}
\begin{equation}\label{quatre}
+\eta(1)(1-\eta(2)) \left[f(\tau_1\eta^{0,1,2})-f(\eta)\right]
\end{equation}
\begin{equation}\label{cinq}
+\eta(-1)(1-\eta(1)) \left[f(\tau_1\eta^{-1,0,1})-f(\eta)\right].
\end{equation}

\noindent Part (\ref{un}) is simple exclusion involving sites away
from the origin, part (\ref{deux}) corresponds to the
``classical'' tagged particle for simple exclusion. Part
(\ref{trois}) is a ``strictly'' 2 steps exclusion involving sites
away from the origin. Part (\ref{quatre}) describes the
``pushing'' of the tagged particle. Finally in part (\ref{cinq})
the tagged particle is ``pushed'' by another particle.

A straightforward adaptation of the simple exclusion case (see
Ferrari (1986)) gives

\begin{theorem}
The Palm measure $\widehat{\nu}_{\alpha}$ of $\nu_{\alpha}$ ({\it
i.e.} the measure on ${\bf X}$ defined by
$\widehat{\nu}_{\alpha}(.)=\nu_{\alpha}(.|\eta(0)=1)$) is
invariant and ergodic for the $k$-step exclusion process as seen
from a tagged pushed particle.
\end{theorem}

\noindent {\bf Sketch of proof:} For instance when $k=2$ it is
enough to show
\[
\mu\{\eta(-1)\eta(0)(1-\eta(1))f(\eta^{-1,0,1})\}=
\mu\{\eta(0)\eta(1)(1-\eta(2))f(\tau_1\eta^{0,1,2})\}
\]
and
\[
\mu\{\eta(-2)\eta(-1)(1-\eta(0))f(\eta^{-2,-1,0})\}=
\mu\{\eta(-1)\eta(0)(1-\eta(1))f(\tau_1\eta^{-1,0,1})\}
\]
which are obvious for any $\mu\in{\cal S}$ (recall that ${\cal S}$
is the set of translation invariant measures on ${\bf X}$).$\carn$

\begin{theorem}
Law of large numbers for the Tagged Pushing Particle ($k=2$).

For a $2$-step exclusion process with initial distribution
$\nu_{\alpha}$, if $Y(t)$ denotes the position at time $t$ of a
Tagged Pushing Particle starting at the origin then
\[
\lim_{t\to\infty}\frac{Y(t)}{t}=(1-\alpha)(1+2\alpha)\ {{\bf
P_{\widehat{\nu}_{\alpha}}}}a.s.
\]
\end{theorem}

\noindent {\bf Proof:} Using the notation of Ferrari (1992a) pp.
41-43 we define the instantaneous increment of the position of the
Tagged Pushing Particle by
\[
\psi(\eta):=\lim_{h\to 0}\frac{{\bf E}(Y(t+h)-Y(t)|Y_t=x)}h
\]
\[
=(1-\eta(x+1))+\eta(x-1)(1-\eta(x+1))+\eta(x+1)(1-\eta(x+2)).
\]
So
\[
\int\psi\ d\nu_{\alpha}=1-\alpha+2\alpha(1-\alpha).
\]
Furthermore $\lim_{h\to 0}{\bf
E}_{\nu_{\alpha}}(Y(t+h)-Y(t))^2/h<+\infty$ so that the conditions
of theorem 9.2 of Ferrari (1992a) are satisfied, which gives the
result.$\carn$

\begin{remark}
\end{remark}
Referring to results of Guiol (1999), we notice that a Tagged
Pushing Particle $X_t$  behaves as a regular tagged particle in
the $k$-step exclusion. Indeed, intuitively, the ``regular''
tagged particle can make long jumps, so is expected to move
faster, but it cannot be pushed; and the rate at which the Tagged
Pushing Particle moves compensates exactly those long jumps.

\section{Generalizations}

\noindent 1. We obtained conservation of local equilibrium for the
Riemann case. It can be generalized to an initial product measure
with a non-constant profile of bounded variation following the
papers of Rezakhanlou (1991) and Landim (1993) (see also Kipnis \&
Landim (1999) chapters 8 and 9). Indeed Condition E is equivalent
to Kru\v{z}kov entropy inequality (see Godlewski \& Raviart(1991),
Lemma 6.1 p.88) which is the key tool of the latter proofs.
Moreover, since we are working with a one-dimensional totally
asymmetric process, we can equivalently consider its
interpretation defined by (\ref{pushing}); this way, particles do
not jump over each other. We can therefore follow Section 6 of
Rezakhanlou (1991) which contains a two-block estimate valid only
for one-dimensional nearest neighbor processes, and enables to
conclude the derivation of hydrodynamic equation without using
Young measures. The key point consists in coupling two versions of
the process, to prove that the number of sign changes between them
can only decrease, due to attractivity and the reduction to a
one-dimensional nearest neighbor case. For $k$-step exclusion,
this result is a straightforward adaptation of Liggett (1976),
Lemma 5.1. For the same reason, extension to deterministic initial
configurations is possible, using Venkatsubramani (1995).\par

\bigskip

\noindent 2. We treated only totally asymmetric $k$-step exclusion
but our proof is also valid for nearest neighbor asymmetric
transition rates, which lead to a flux function with the same
properties. However one has to compute carefully this function,
due to the complexity of the $k$-step dynamics. For instance for
$k=5$ and transition rates $p(x,x+1)=p$, $p(x,x-1)=q=1-p$, $p>q$,
we obtain
\[
\begin{array}{ll}
G_5^{p,q}(u)=&u(1-u)\left[(p-q)\left\{(1+2u)+3u^2(1-pq)\right.\right.\\
&\left.\left.+4u^3(1-2pq)+5u^4(1-3pq+p^2q^2)\right\}+3u^2
2p^4q\right].\\
\end{array}
\]
Indeed the last term corresponds to a ``cycle'' jump: To go from
$x$ to $x+3$ when sites $x+1$ and $x+2$ are occupied, the particle
follows the path $(x,x+1,x+2,x+1,x+2,x+3)$.\par The extension to
an initial product measure with a general profile is still
possible, using Rezakhanlou (1991) and Landim (1993), but with the
help of Young measures.\par

\bigskip

\noindent {\bf Aknowledgments:} We thank E. Andjel, C. Bahadoran,
P. Ferrari, J. Krug and G. Sch\"utz for fruitful discussions. This
work was supported by the agreement USP/COFECUB n$^o$ UC 45/97,
and FAPESP project n$^o$ 1999/04874-6. The three authors thank IME
at S\~ao Paulo for their kind hospitality, H.G. \& K.R. thank
Universit\'e de Rouen where part of this work was realized.

\end{document}